\input epsf.sty
\documentstyle{amsart} 
\newcommand{\bd}{\text{Bd}\,}
\newcommand{\be}{\begin{enumerate}}

\newcommand{\ca}{{\cal A}}
\newcommand{\cb}{{\cal B}}
\newcommand{\cd}{{\cal D}}
\newcommand{\cf}{{\cal F}}
\newcommand{\cl}{{\cal L}}
\newcommand{\cm}{{\cal M}}
\newcommand{\cn}{{\cal N}}
\newcommand{\cp}{{\cal P}}
\newcommand{\cs}{{\cal S}}
\newcommand{\cu}{{\cal U}}
\newcommand{\cv}{{\cal V}}
\newcommand{\czz}{{\cal Z}{\cal Z}}
\newcommand{\da}{\Delta}
\newcommand{\ee}{\end{enumerate}}

\newcommand{\hocolim}{\text{\bf hocolim}\,}
\newcommand{\hra}{\hookrightarrow}

\newcommand{\lw}{\overleftarrow{W$\,\,$}}

\newcommand{\piv}{\text{piv}\,}
\newcommand{\pr}{{\bf Proof. }}
\newcommand{\ra}{\rightarrow}
\newcommand{\rn}{{\Bbb R}^n}

\newcommand{\rw}{\overrightarrow{W$\,\,$}}
\newcommand{\sm}{\setminus}

\newcommand{\topo}{\text{\bf Top}}

\newtheorem{thm}{Theorem}[section]
\newtheorem{df}  [thm]{Definition}

\newtheorem{lm}  [thm]{Lemma}

\newtheorem{prop}[thm]{Proposition}

\numberwithin{equation}{section}

\begin{document}

\title[A comparison of Vassiliev and Ziegler-\v{Z}ivaljevi\'{c} models]
{A comparison of Vassiliev and Ziegler-\v{Z}ivaljevi\'{c} models
for homotopy types of subspace arrangements}
              \author{Dmitry N. Kozlov}
  \date{\today. \\[0.05cm]
 \hskip15pt This research was supported by the Research Grant 
dnr~2247/1999 of the Swedish Natural Science Research Council.
}

\address{Department of Mathematics, Royal Institute of Technology, 
Stockholm, S-100 44, Sweden. 
\newline
{\it current address:} Department of Mathematics, University of Bern, 
Sidlerstrasse 5, Bern, CH-3012, Switzerland.}

\email{kozlov@@math.kth.se, kozlov@@math.ethz.ch}

\begin{abstract} 
  In this paper we represent the Vassiliev model for the homotopy type 
of the one-point compactification of subspace arrangements 
as a homotopy colimit of an appropriate diagram over the nerve complex
of the intersection semilattice of the arrangement. Furthermore, using 
a~generalization of simplicial collapses to diagrams of topological 
spaces over simplicial complexes, we construct an~explicit deformation 
retraction from the Vassiliev model to the Ziegler-\v{Z}ivaljevi\'{c} 
model.
\end{abstract}

\maketitle

                 \section{Introduction}

  Goresky and MacPherson, \cite[Part III]{GM}, were the first 
to express the cohomology groups of the complement of a subspace
arrangement~$\ca$ in terms of the homology groups of the order 
complexes of lower intervals of the associated intersection 
semilattice. Following that, there was a sizable body of work 
studying the topological properties of the complement of subspace 
arrangements, or, dually, of the one-point compactification of 
the union of subspaces, which we denote by $\widehat\cu(\ca)$, 
see~\cite{Bj1,Hu,V,Z,ZZ}. Especially elucidating argument can 
be found in \cite[Chapter~II.5]{Z}.

  In particular, two models were constructed, one by Vassiliev, 
\cite{V}, and one by Ziegler and \v{Z}ivaljevi\'{c}, \cite{ZZ}, 
reproducing $\widehat\cu(\ca)$ up to homotopy equivalence. 
The Ziegler-\v{Z}ivaljevi\'{c} model is based on the notion of 
homotopy colimit, dating back at~least to \cite{BK}, but see also
\cite{WZZ} for a fresh approach; while Vassiliev's construction
is explicitly geometrical. It was explicitly verified 
in~\cite[page 140]{WZZ} that the two models are homotopy equivalent.

  The purpose of this paper is twofold. First we find a~presentation
for the Vassiliev model as a~certain homotopy colimit, thus bringing
the two models to a~common formal framework. Second, by using
a~diagram-theoretic generalization of simplicial collapses,
coupled with the technical machinery of Discrete Morse Theory, 
\cite{Fo}, we describe a~sequence of generalized collapses leading 
from the Vassiliev model to the Ziegler-\v{Z}ivaljevi\'{c} model.
This, in turn, connects the two models by a deformation retraction.

\vskip4pt
\noindent
{\bf Acknowledgments.} I am grateful to Eva-Maria Feichtner for the 
careful proofreading of this paper, and to G\"unter M.~Ziegler for 
the useful comments on the history of the subject.


         \section{Background}
    
  \subsection{The terminology of posets}   $\,$
\vskip3pt

A poset is a set with a specified partial order. We say that a poset
$P$ is a {\bf semilattice} if for any $x,y\in P$ the sets 
$\{z\in P\,|\,x\leq z,y\leq z\}$, resp.~$\{z\in P\,|\,x\geq z,y\geq z\}$
are either empty or have minimal, resp.~maximal elements.

  Let $\cp$ denote the full subcategory of the category of all small 
categories consisting of posets. Here posets are viewed as categories
in the standard way, i.e., with elements being the objects
and order relations being the morphisms. Let furthermore $\topo$
denote the category of topological spaces and continuous maps.
 
The definition of the nerve of a category goes back to Quillen, \cite{Q},
and Segal, \cite{Se}, we state it only in the special case of posets,
and we also compose it at once with the functor mapping simplicial 
complexes to their geometric realizations.

\begin{df}
  The functor $\da:\cp\ra\topo$ maps a poset $P$ to the geometric 
realization of the simplicial complex whose vertices are 
the elements of $P$ and whose simplices correspond to chains 
(totally ordered subsets) of $P$. $\da(P)$ is commonly 
known as the {\bf order complex} of $P$.
\end{df}
 
  For $x\in P$, we denote by $P_{\leq x}$ the full subposet of $P$
consisting of elements $\{y\in P\,|\,y\leq x\}$. Analogously,
$P_{<x}$ is the full subposet of $P$ consisting of elements
$\{y\in P\,|\,y<x\}$.

  The {\bf barycentric subdivision} of a poset $P$, denoted $\bd(P)$ 
is a poset whose elements are all non-empty chains of $P$ partially 
ordered by inclusion.

  Given a simplicial complex $K$, we denote by $\cf(K)$ its 
{\bf face poset}, which is the poset consisting of all non-empty 
faces of $K$ partially ordered by inclusion.

For $x,y\in P$, $x\geq y$, we denote by $I(y\hra x)$ the inclusion map
of the simplicial complexes $I(y\hra x):\da(P_{\leq y})\hookrightarrow
\da(P_{\leq x})$. 

  \subsection{The terminology of subspace arrangements}   $\,$
\vskip3pt
 
 A {\bf subspace arrangement} is a collection $\ca=\{A_1,\dots,A_k\}$ 
of affine linear subspaces in $\rn$, such that if $A_i\subseteq A_j$,
then $A_i=A_j$. To this collection we associate the following invariants:
\begin{itemize}

\item The {\bf intersection semilattice} $\cl(\ca)$ consisting of all
possible non-empty intersections of $A_i$'s ordered by reverse inclusion;

\item The collection $\cb(\ca)=\{B(x)\,|\,x\in\cl(\ca)\}$ of corresponding
affine subspaces indexed by the elements of the intersection semilattice;

\item We denote $\cu(\ca)=\cup_{i=1}^k A_i$ and $\cm(\ca)=\rn\sm\cu(\ca)$. 
Let $\widehat\cu(\ca)$ denote the one-point compactification of $\cu(\ca)$. 
\end{itemize}

In the rest of this section, following Vassiliev and 
Ziegler-\v{Z}ivaljevi\'{c}, \cite{V,ZZ}, we define two different 
topological spaces both of which are homotopy equivalent to 
$\widehat\cu(\ca)$ (in particular, they are of course homotopy equivalent 
to each other).
   
\subsection{Homotopy colimits} \label{ss2.3}  $\,$
\vskip3pt

\begin{df} A {\bf diagram of topological spaces over a poset $P$},  
is a~covariant functor from $P$ to $\topo$.
\end{df}

If the functor is denoted by $\cd$, and $x$ is an element of $P$, 
we use $\cd(x)$ to denote the topological space associated to $x$;
and if $x,y\in P$, $x\geq y$, we use $\cd(x\ra y)$ to denote the continuous 
map associated to the order relation $x\geq y$ (which is a morphism in $P$
viewed as a category).

In this paper the topological spaces $\cd(x)$ are always direct products
of (geometric realizations of) simplicial complexes with linear subspaces, 
and the maps $\cd(x\ra y)$ are always inclusions.

\begin{df} \label{hocolimdf}
  The {\bf homotopy colimit} of a diagram of topological spaces 
$\cd:P\ra\topo$, denoted by $\hocolim(\cd)$, is the colimit of 
the functor $\da(\cd):\bd(P)\ra\topo$ defined by: 
\begin{itemize}
\item on the elements: 
   $\da(\cd)(x_1>\dots>x_t)=\da(P_{\leq x_t})\times\cd(x_1)$;
\item on the morphisms: 

$\da(\cd)((x_1>\dots>x_t)\ra(x_{i_1}>\dots>x_{i_p}))=
I(x_t\hra x_{i_p})\times\cd(x_1\ra x_{i_1})$.
\end{itemize}
\end{df}

  One of the main sources for details on homotopy colimits
is \cite{BK}, see also \cite{WZZ} for many combinatorial applications
of the concept.
 
 \vskip4pt
  Later on, we shall need the following explicit description of 
the topological space $\hocolim(\cd)$. Consider the disjoint
union of spaces $\cd(x)$, for $x\in P$, then for any order relation
$x>y$ glue in the mapping cylinder of the map $\cd(x\ra y)$, taking
$\cd(x)$ as the source, and $\cd(y)$ as the base of it;
for every triple $x>y>z$ glue in the "mapping triangle"
of maps $\cd(x\ra y)$ and $\cd(y\ra z)$ and so on through the entire
order complex of $P$. Of course, while geometrically intuitive, this 
description follows word-by-word the definition of the colimit.
\vskip4pt 
 
  An important special example which we need in this paper is 
the case when $P$ is the face poset of a simplicial complex $K$, 
$P=\cf(K)$. In this case, we call $\cd:P\ra\topo$, a~{\bf diagram 
over the simplicial complex $K$.} 

\begin{df} \label{df2.4}
  Let $\cd:P\ra\topo$ be a diagram of topological spaces over 
a poset, define a diagram over the simplicial complex $\da(P)$,
$\bd(\cd):\bd(P)\ra\topo$ as follows:
\begin{itemize}
\item on objects: $\bd(\cd)(x_1>\dots>x_k)=\cd(x_1)$;
\item on morphisms: $\bd(\cd)((x_1>\dots>x_k)\ra 
(x_{i_1}>\dots>x_{i_t}))=\cd(x_1\ra x_{i_1})$.
\end{itemize}
\end{df}

  As the next proposition shows (verification is left to the reader) 
any diagram over a poset can be replaced with a diagram over 
a~simplicial complex.

\begin{prop}\label{prop2.5}
   For any diagram $\cd$ of topological spaces over a~poset,
 the space $\hocolim(\bd(\cd))$ is homeomorphic to $\hocolim(\cd)$.
\end{prop}  
 
 \section{Description of the models. 
Representing the Vassiliev model \\ as a~homotopy colimit}
 
\subsection{Ziegler-\v{Z}ivaljevi\'{c} model} \label{ss2.1} $\,$
\vskip3pt

The following diagram was suggested for consideration in \cite{ZZ,Z}.

\begin{df}
  Given an affine subspace arrangement $\ca$ in $\rn$, 
the diagram $\czz(\ca):\cl(\ca)\ra\topo$ is defined by:
\begin{itemize}
\item on objects: $\czz(\ca)(x)=B(x)$;
\item on morphisms: $\czz(\ca)(x\ra y)$ is the corresponding inclusion 
map of $B(x)$ into $B(y)$.
\end{itemize}
\end{df}

  It follows from Proposition~\ref{prop2.5} that 
$\hocolim(\czz(\ca))$ is homeomorphic to the homotopy colimit of 
the corresponding diagram over the simplicial complex $\da(\cl(\ca))$.
The following proposition is a consequence of the Projection Lemma,
\cite[XII.3.1(iv)]{BK}, see~\cite{ZZ,Z}.
\begin{prop}
   For an affine subspace arrangement $\ca$ in $\rn$, $\widehat\cu$
is homotopy equivalent to $\hocolim(\czz(\ca))\cup\{\infty\}$ . 
\end{prop}
  
  By using the Homotopy Lemma, \cite[XII.4.2]{BK}, Ziegler and 
\v{Z}ivaljevi\'{c} could then prove the following formula for 
the homotopy type of~$\widehat\cu(\ca)$.
  
\begin{thm}
  For an affine subspace arrangement $\ca$ in $\rn$
  $$\widehat\cu(\ca)\simeq\vee_{x\in\cl(\ca)}
     (\da(\cl(\ca)_{<x})*S^{\dim(B(x))}).$$
\end{thm}  
  
  And hence, by Alexander duality, one gets the cohomology groups of 
the complement, originally due to Goresky and MacPherson, \cite{GM}.
  
\begin{thm}
  For an affine subspace arrangement $\ca$ in $\rn$
  $$\widetilde H^i(\cm(\ca);{\Bbb Z})\cong\oplus_{x\in\cl(\ca)}
  \widetilde H_{n-i-\dim(B(x))-2}(\da(\cl(\ca)_{<x};{\Bbb Z}).$$
\end{thm}

 \subsection{Vassiliev model}$\,$
\vskip3pt

 Vassiliev has suggested a slightly different modification of
the subspace arrangement. The idea is to 
"simplicially blow up" the intersections of the subspaces.
Vassiliev calls it a~{\it geometric resolution}.

More precisely: take $N$ to be a sufficiently large number and
embed subspaces $A_i$ into ${\Bbb R}^N$ in a generic position;
for every $x\in\cu(\ca)$, let $V(x)$ be the convex hull
of the images of $x$ in ${\Bbb R}^N$. 

Let $V(\ca)=\cup_{x\in\cu(\ca)}V(x)$. It is a~"resolution" of 
the arrangement in the following sense. 

\begin{lm} \cite[Lemma 1, p.~120]{V}
  One can choose $N$ sufficiently large, and the embedding 
sufficiently generic, so that, for every $x\in\cu(\ca)$, $V(x)$ is 
a~simplex with vertices being the images of $x$ in ${\Bbb R}^N$, 
and, for every $x,y\in\cu(\ca)$, $x\neq y$, the simplices $V(x)$ and 
$V(y)$ do not intersect.
\end{lm}
 
\begin{prop} \cite[Lemma 2, p.~120]{V}
   The one-point compactification of the geometric resolution
$V(\ca)\cup\{\infty\}=\widehat V(\ca)$ is homotopy equivalent to 
$\widehat\cu(\ca)$.
\end{prop} 

  Vassiliev then, by means of an explicit argument using Stratified 
Morse Theory of Goresky and MacPherson, \cite{GM}, obtains 
a~description for the homotopy type of $\widehat\cu(\ca)$
which is essentially identical to the Ziegler-\v{Z}ivaljevi\'{c}
description. Amazingly both results were obtained simultaneously 
and independently.

  An observation which both Vassiliev and Ziegler-\v{Z}ivaljevi\'{c}
make is that it follows by Spanier-Whitehead duality that the stable 
homotopy type of $\cm(\ca)$ is defined by the combinatorial data of 
the arrangement (the intersection semilattice together with 
the dimension information), while is it well-known that the homotopy 
type of $\cm(\ca)$ is not a combinatorial invariant, see~\cite{V,Z}, 
\cite[Theorem 3.4]{ZZ}.

\subsection{Representing Vassiliev model as a homotopy colimit}$\,$
\vskip3pt

\begin{df}
  Given a semilattice $P$, we define the simplicial complex
$\cn(P)$ as follows:
\begin{itemize}
\item the vertices of $\cn(P)$ are the minimal elements of $P$;
\item the simplices of $\cn(P)$ are those collections of minimal 
elements of $P$ which have a join in $P$.
\end{itemize}
$\cn(P)$ is known as the {\bf nerve complex} of $P$.
\end{df}

  It was proved by Leray, \cite{Le}, that the C{\^e}ch homology
groups of $\cn(P)$ and of $\da(P)$ are equal, and by Borsuk,
\cite{Bo}, that the two complexes are actually homotopy equivalent. 

  Next, we use the notion of the~nerve complex of the intersection
lattice to define a~specific diagram of spaces associated to 
an~affine subspace arrangement, which to our knowledge 
was not previously considered in the literature. 

\begin{df}
  Let $\ca=\{A_1,\dots, A_k\}$ be an affine subspace arrangement, and 
denote the elements of $\cl(\ca)$ corresponding to $A_1,\dots, A_k$
by $a_1,\dots,a_k$. We define the~{\bf Vassiliev diagram} 
$\cv:\cf(\cn(\cl(\ca)))\ra\topo$ to be the functor specified by:
\begin{itemize}
\item on elements: $\cv(\{a_{i_1},\dots,a_{i_k}\})=
A_{i_1}\cap\dots\cap A_{i_k}$;
\item on morphisms: the maps are inclusions 
$$\cv(\{a_{i_1},\dots,a_{i_k}\}\ra\{a_{j_1},\dots,a_{j_q}\})=
(A_{i_1}\cap\dots\cap A_{i_k})\hookrightarrow(A_{j_1}\cap\dots\cap A_{j_q}),$$
for any $\{j_1,\dots,j_q\}\subseteq\{i_1,\dots,i_k\}$.
\end{itemize}
\end{df}

\begin{thm}  \label{main1}
  $\hocolim(\cv)\cup\{\infty\}$ is homeomorphic to $\widehat V(\ca)$.
\end{thm}
\pr
  It is immediate from the definitions that $\hocolim(\cv)\cup\{\infty\}$
is a~"bary\-centric subdivision" of $\widehat V(\ca)$, that is, all
the simplices which Vassiliev spans on the images of points under
the generic embedding are barycentrically subdivided in 
$\hocolim(\cv)\cup\{\infty\}$. Other than that, there is no difference in 
the construction and so we conclude that the two spaces are homeomorphic.
\qed
\vskip4pt

  Again, it follows from Proposition~\ref{prop2.5}, that 
$\hocolim(\cv)$ is homeomorphic to the homotopy colimit of the 
corresponding diagram over the simplicial complex $\bd(\cn(\cl(\ca)))$.

  \section{A deformation retract from the Vassiliev model to 
the Ziegler-\v{Z}ivaljevi\'{c} model}

 \subsection{Single collapse}$\,$
\vskip3pt

  Assume that we have a diagram over a simplicial complex $K$,
$\cd:\cf(K)\ra\topo$, such that for some simplices 
$\sigma,\tau\in\cf(K)$ the following is true:
\begin{itemize}
\item $\sigma<\tau$, and there exists no simplex in $K$, other than
$\tau$ and $\sigma$ itself, which contains $\sigma$, in particular $\tau$ 
is maximal; in such situation one says that removing $\sigma$ and $\tau$
from $K$ is an elementary simplicial collapse;
\item $\cd(\tau\ra\sigma)$ is an identity map.
\end{itemize}

\begin{prop} \label{prop4.1}
  In the situation above there exists a deformation retract from 
$\hocolim\cd$ to $\hocolim\cd'$, where $\cd':\cf(K\sm\{\sigma,\tau\})
\ra\topo$ is the restriction of the functor $\cd$.
\end{prop}
\pr 
  The desired retract is a simple generalization of the deformation
which retracts a mapping cylinder to the target space. It can be
easily visualized as follows: think that we have a string connecting
the unique vertex $v$ of $\tau$ which does not lie in $\sigma$
to the barycenter $w$ of $\sigma$, and that we start to shrink the 
string so that $w$ approaches $v$ over an interval of time $[0,1]$
($w$ coincides with $v$ at moment 1). We let the entire homotopy
colimit be deformed accordingly, and refer to the explicit 
description of homotopy colimits in Subsection~\ref{ss2.3} 
for visualizing this process. 

  This is clearly a retract from $\hocolim\cd$ to $\hocolim\cd'$. 
The continuity of this deformation at any time $0\leq t<1$ follows 
from the fact that $\cd(\tau\ra\sigma)$ is an identity map, and 
the continuity at $t=1$ follows from the definition of the category 
$\topo$ (the morphisms are continuous maps).
\qed

 \subsection{Terminology of Discrete Morse Theory.}$\,$
\vskip3pt

  Although unaware of an exact reference, we are confident that
it is folklore knowledge that for every finite semilattice $P$
there is a sequence of collapses leading from $\bd(\cn(P))$ to
$\da(P)$. However, to use Proposition~\ref{prop4.1},
we need to check a condition that certain maps are identities, 
so we will list this sequence of collapses explicitly.

  It is handy to use the formal setup of Discrete Morse Theory.
We provide below the necessary terminology and results
for the special case that we need, see~\cite{Fo} for further details.

Let $K$ be a~simplicial complex. A~{\bf matching} $W$ on $P=\cf(K)$ 
(cf.~\cite[Definition~9.1]{Fo}) is a~set of disjoint pairs 
$(\sigma,\tau)$ such that $\tau,\sigma\in P$, $\tau\succ\sigma$, 
(``$\succ$'' denotes the covering relation). We set 
$$\rw=\{\sigma\in P\,|\, \text{ there exists } \tau \text{ such that } 
(\sigma,\tau)\in W\},$$ 
$$\lw=\{\tau\in P\,|\, \text{ there exists } \sigma \text{ such that } 
(\sigma,\tau)\in W\}.$$ 
If $(\sigma,\tau)\in W$ then we set $W(\sigma)=\tau$.

\begin{df} (cf. \cite[Definition~9.2]{Fo}).
  A~matching is called {\bf acyclic} if it is impossible to find a~sequence 
$\sigma_0,\dots,\sigma_t\in\rw$, such that $\sigma_0\neq \sigma_1$, 
$\sigma_0=\sigma_t$, and $W(\sigma_i)\succ\sigma_{i+1}$, for $0\leq i\leq t-1$.
\end{df}

  The following proposition is the only fact that we need for our 
argument, see also \cite[Corollary~3.5,Theorem 9.3]{Fo}, and 
\cite[Theorem 3.2 (2)]{Ko}.

\begin{prop}\label{morse}
  Let $K$ be a~simplicial complex and $P=\cf(K)$ be its face poset. 
Let $W$ be an acyclic matching on $P$. If the unmatched simplices 
form a~subcomplex $K^C$ of $K$, then there is a~sequence of elementary 
collapses leading from $K$ to $K^C$.
\end{prop}

 \subsection{An acyclic matching for our case.}$\,$ \label{ss4.3}
\vskip3pt

  Let $P$ be a semilattice. We call a~set $\{a_1,\dots,a_t\}
\subseteq\min(P)$ {\it complete} if 
\begin{itemize}
\item $\vee_{i=1}^t a_i$ exists;
\item if $x\leq \vee_{i=1}^t a_i$, and $x\in\min(P)$, then 
$x\in\{a_1,\dots,a_t\}$; in other words 
$\min(P)\cap(P_{\leq\vee_{i=1}^t a_i})=\{a_1,\dots,a_t\}$.
\end{itemize}

  Otherwise a subset of $\min(P)$ is called {\it incomplete}.
For any subset $\{b_1,\dots,b_q\}\subseteq\min(P)$, such that
$\vee_{j=1}^q b_j$ exists, we call $\min(P)\cap(P_{\leq\vee_{j=1}^q b_j})$ 
the~{\it completion} of $\{b_1,\dots,b_q\}$, and denote it by 
$C(\{b_1,\dots,b_q\})$. Clearly, a set is complete iff it is equal
to its own completion.

  By construction, $\da(P)$ is the full subcomplex of $\bd(\cn(P))$
spanned by the vertices which are enumerated by the complete subsets
of $\min(P)$. 

  Let us now define an acyclic matching on $\bd(\cn(P))$.
For a simplex $\Sigma=(S_1<\dots<S_t)$ of $\bd(\cn(P))$ let
$\piv(\Sigma)$ denote the incomplete set $S_i$ with the maximal possible 
index $i$, if it exists; set $\piv(\Sigma)=\emptyset$ if it does not.
If $\piv(\Sigma)\neq\emptyset$, set $\iota(\Sigma)$ to be equal to 
the index of $\piv(\Sigma)$ in $\Sigma$. Define
$$\rw=\{\Sigma=(S_1<\dots<S_t)\,|\,\piv(\Sigma)\neq\emptyset
\text{ and }C(\piv(\Sigma))\not\in\Sigma\}.$$
Correspondingly we define
$$\lw=\{\Sigma=(S_1<\dots<S_t)\,|\,\piv(\Sigma)\neq\emptyset
\text{ and }C(\piv(\Sigma))\in\Sigma\}.$$
Finally, for $\Sigma\in\rw$ we define 
$W(\Sigma)=\Sigma\cup \{C(\piv(\Sigma))\}$. 
Clearly $\bd(\cn(P))=\da(P)\cup\rw\cup\lw$ and the union is disjoint.
  
\begin{prop}
  The matching $W$ described above is acyclic.
\end{prop}
\pr
  Assume that there exists a~sequence $\Sigma_0,\dots,\Sigma_t\in\rw$, 
such that $\Sigma_0\neq \Sigma_1$, $\Sigma_0=\Sigma_t$, and 
$W(\Sigma_i)\succ\Sigma_{i+1}$, for $0\leq i\leq t-1$. 
We have the following equalities and inequalities: 
$$\iota(\Sigma_0)=\iota(W(\Sigma_0))>\iota(\Sigma_1)=\iota(W(\Sigma_1))>
\dots =\iota(W(\Sigma_{t-1}))>\iota(\Sigma_t)=\iota(\Sigma_0),$$
which yields a contradiction.
\qed

 \subsection{The deformation retract theorem}$\,$
\vskip3pt

\begin{thm} \label{main2}
  $\hocolim(\czz)\cup\{\infty\}$ is a deformation retract of 
  $\hocolim(\cv)\cup\{\infty\}$. 
\end{thm}
\pr
  It is enough to show that $\hocolim(\bd(\czz))$ is a deformation 
retract of $\hocolim(\bd(\cv))$. For that we need to verify that
in the matching described in Subsection~\ref{ss4.3} the maps
within the matched pairs are always identities.

  Since both diagrams are obtained by subdivisions, it follows
from Definition~\ref{df2.4} that the desired maps are
obviously identities in all cases, except possibly when a~pair 
$(\Sigma,W(\Sigma))$ is such that $\piv(\Sigma)$ is the maximal
element of $\Sigma$. 

  In this case, if we use the notations $\piv(\Sigma)=\{a_1,\dots,a_t\}$, 
and $C(\piv(\Sigma))=\{a_1,\dots,a_t,a_{t+1},\dots,a_{t+k}\}$, then
the desired map is the inclusion $\cap_{i=1}^{t+k}A_i\hookrightarrow
\cap_{i=1}^t A_i$, which is the identity by definition of the completion 
(here, $A_i\in\ca$ denotes the subspace indexed by $a_i\in\cl(\ca)$).
\qed

 The deformation procedure is illustrated on Figure 1 for the
example of the arrangement consisting of 3 lines, all intersecting
in the same point.

$$\epsffile{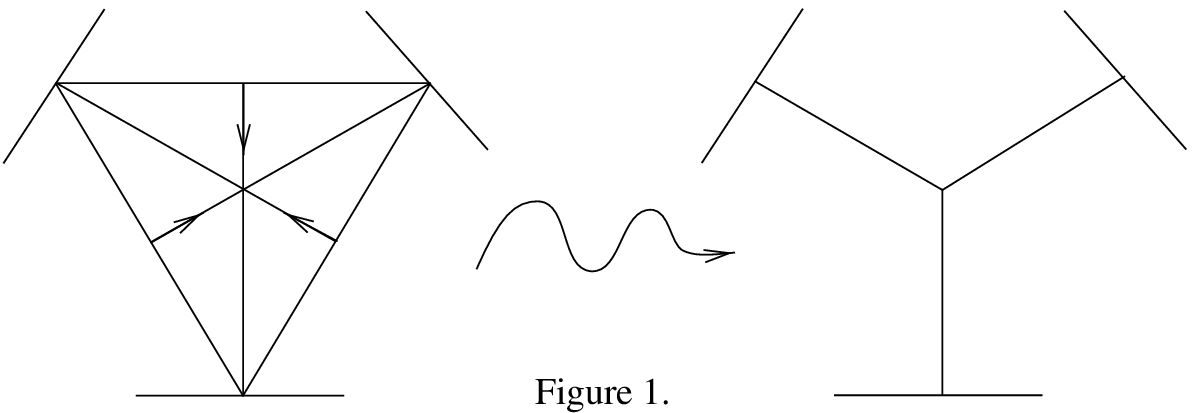}$$

  \subsection{Final remark}$\,$
  \vskip3pt
  
  Removing the infinity throughout the paper yields the uncompactified
version of the result.

\end{document}